\documentclass[a4paper,11pt]{amsart}
\usepackage{amssymb,amsfonts,amsxtra, mathrsfs,placeins,graphicx,verbatim}
\usepackage[all]{xy}
\xyoption{line}
\usepackage{fullpage}
\usepackage{euscript}
\newtheorem{theorem}{Theorem}[section]
\newtheorem{cor}[theorem]{Corollary}
\newtheorem{lem}[theorem]{Lemma}
\newtheorem{prop}[theorem]{Proposition}

\theoremstyle{definition}

\newtheorem{example}[theorem]{Example}
\newtheorem{defi}[theorem]{Definition}
\newtheorem{rem}[theorem]{Remark}

\numberwithin{equation}{section}

\DeclareMathOperator{\Hom}{Hom}

\DeclareMathOperator{\im}{Im}
\DeclareMathOperator{\Ker}{Ker}

\DeclareMathOperator{\Tr}{Tr}

\def\E{\mathcal E}
\def\P{\mathcal P}
\def\B{\mathsf B}

\def\End{\operatorname{End}}

\def\ra{\rightarrow}

\def\S{\mathbb{S}}

\def\ground{\mathbf k}

\newcommand{\noproof}{\begin{flushright}\ensuremath{\square}\end{flushright}}

\newcommand{\BV}{\mathsf{BV}}

\newcommand{\bv}{\mathsf{bv}}

\def\O{\mathcal O}
\def\E{\mathcal E}

\def\ev{\operatorname{ev}}

\def\id{\operatorname{id}}

\thanks{}

\begin{document}

\title[Minimal models...]{Abstract Hodge decomposition and minimal models for cyclic algebras}
\author{J. Chuang  \and A.~Lazarev}
\thanks{The first author is supported by an EPSRC advanced research fellowship.
The second author is partially supported by an EPSRC research grant}
\address{Centre for Mathematical Science\\City University\\London EC1V 0HB\\UK}
\email{J.Chuang@city.ac.uk}
\address{University of Leicester\\ Department of Mathematics\\Leicester LE1 7RH, UK.}
\email{al179@leicester.ac.uk}

\begin{abstract}
We show that an algebra over a cyclic operad supplied with an additional linear algebra datum called \emph{Hodge decomposition} admits a minimal model whose structure maps
are given in terms of summation over trees. This minimal model is unique up to homotopy.
\end{abstract}

\maketitle
\section{Introduction}

Operadic algebras (such as $A_\infty$- , $C_\infty$- and $L_\infty$-algebras) were originally introduced for the needs of homotopy theory but now figure prominently also in algebraic geometry and certain parts of theoretical physics, particularly in those aspects which concern mirror symmetry. One of the highlights of the theory of operadic algebras is the theorem due to Kadeishvili \cite{Kad} on the existence of minimal models for $A_\infty$-algebras.

Relatively recently this theorem was revisited by Merkulov, Markl and Kontsevich-Soibelman \cite{Mer, markl, KS}. They gave explicit formulas for structure maps of minimal models in terms of summations over trees.

In \cite{CL2} we introduced a new approach to the construction of explicit minimal models, giving a conceptual explanation for the Merkulov tree sum formulas. We worked in an operadic setting, so our results apply to homotopy algebras interpreted in the broadest sense, including as special cases $C_\infty$-algebras and $L_\infty$-algebras. Moreover, these methods extended naturally to modular operads, and in this setting minimal models are given by sums over general graphs similar to those appearing in the Feynman diagram expansions of path integrals. As a special case we recovered the theory of minimal models for homotopy algebras equipped with non-degenerate bilinear forms.

The purpose of the present article is to demonstrate how our approach can be adapted for homotopy algebras with bilinear forms which are not necessarily non-degenerate. This is important because the most interesting examples
such as de Rham and Dolbeault algebras are infinite dimensional and thus cannot support a non-degenerate inner product. However their minimal models tend to be finite dimensional and so it is natural to ask whether their inner products could be made compatible with the higher multiplications.

Finite-dimensional models of operadic algebras with a compatible scalar product are also important because of their connection to various moduli spaces; for example Gromov-Witten invariants encode the action of the homology operad of moduli spaces of curves on the Hodge cohomology of smooth projective varieties.

Another source of interest in finite-dimensional models is the the so-called `direct construction' of Kontsevich \cite{K1,HL} which associates to an $A_\infty$-algebra with a non-degenerate scalar product a family of cohomology classes on moduli spaces of Riemann surfaces.

The appropriate setting in which to study these constructions is that of algebras over \emph{cyclic} operads. Our main result states that an algebra $V$ (not necessarily finite dimensional) over a cyclic operad admits an explicit minimal model whose structure maps are given by a Merkulov-type formula \emph{provided the underlying complex of $V$ possesses a Hodge decomposition}. The Hodge decomposition of a complex with an inner product is a certain linear algebra datum which is interesting in its own right and we discuss this notion separately. It is an abstract version of the structure possessed by the De Rham algebra of a smooth manifold or a Dolbeault algebra of a complex manifold. Without the presence of a form this structure is known under the name of \emph{strong deformation retract} data. If the space $V$ is finite dimensional then a Hodge decomposition always exists but in infinite-dimensional situations its existence is not straightforward.

There are several papers which treated the $A_\infty$-case using more ad hoc methods, notably \cite{Pol, Laz, Kajiura}.
Other methods based on noncommutative geometry \cite{KS1,Bal} potentially generalize to the operadic context but do not produce an explicit formulas.

The paper is organized as follows. In the remaining part of the introduction we outline the explicit formula for a minimal model in the $A_\infty$-case. The reason for this is that although our proof does not simplify for this particular choice of an operad, the notational complications make the formulation of the general case considerably harder to digest. Chapter 2 introduces and studies the notion of abstract Hodge decomposition. In Chapter 3 we discuss various notions of a cyclic algebra over a cyclic operad; all these notions coincide in the case of a non-degenerate form. Chapter 4 contains our main result and various examples.

We thank Pavel Hajek for spotting a mistake in Example 2.9 in an earlier version of this paper.

\subsection{Notation and conventions} We work in the category of $\mathbb Z/2$-graded differential graded (dg) vector spaces over a fixed ground field $\ground $ of characteristic zero. For a dg vector space $V$ we denote by $\Pi V$ the parity reversion of $V$; thus $(\Pi V)_0=V_1$ and $(\Pi V)_1=V_0$. The homology of a dg vector space $V$ will be denoted by $H(V)$. The group of permutations of $n$ symbols will be denoted by $S_n$. An $\S$-module is a collection $E(n)$ of dg vector spaces together with a right action of $S_n$ on each $E(n)$. Similarly an  $\S_+$-module is a collection $E(n)$ of dg vector spaces together with a right action of $S_{n+1}$ on each $E(n)$.  A dg operad is an $\S$-module $\O=\{\O(n)\}, n=0, 1,\ldots$ together with right actions of $S_n$ on each $\O(n)$ and structure maps $\circ_i:\O(m)\otimes\O(n)\rightarrow\O(m+n-1)$, where $i=1,\ldots, m$ subject to natural axioms, see, e.g. \cite{MSS}. We will denote by $\phi^{\sigma}$ the right action of a permutation $\sigma$ on $\phi$ and by $z(n)$ the $n+1$-cycle $z(n)=(0,1,\ldots,n)\in S_{n+1}$.

\subsection{Minimal models for cyclic $A_\infty$-algebras}
We will now explain how our construction works for $A_\infty$-algebras.
Recall that an $A_\infty$-algebra is a dg vector space $V$ together with odd maps
$$m_n:(\Pi V)^{\otimes n}\to \Pi V,   \qquad n\geq 2,$$
 such that
\begin{eqnarray}\label{eq:ainf}
\sum
_{i+j+k=n}
m_{i+1+k} (\id^{\otimes i}\otimes m_j \otimes \id^{\otimes k})=0, \qquad n\geq 1,
\end{eqnarray}
where $m_1$ is the differential of $V$.

Suppose $V$ is equipped with a symmetric  bilinear form (even or odd)
$\langle,\rangle : V\otimes V \to \ground$.
Then $(V,\{m_i\})$ is called a \emph{cyclic} (or \emph{symplectic}) $A_\infty$-algebra
if the tensor $\tilde{m}_n\in((\Pi V)^{\otimes n})^*$ defined by the formula
\begin{equation}\label{eq:cyclic}
\tilde{m}_n(\Pi v_0,\ldots, \Pi v_n)
=\langle m_n(\Pi v_0,\ldots, \Pi v_{n-1}), \Pi v_n \rangle
\end{equation}
is invariant with respect to all cyclic permutations on $((\Pi V)^{\otimes n})^*$.

Now let $(V,\{m_i\})$ be a cyclic $A_\infty$-algebra. We explain how a cyclic $A_\infty$-algebra
structure is induced on the homology
$H(V)$, a la Merkulov/Kontsevich-Soibelman. To implement the construction we choose a decomposition
$V=X \oplus \im(d) \oplus W$ such that $X\oplus \im(d) = \ker(d)$ and $W$ is orthogonal
to $X\oplus W$. When $V$ is finite dimensional and also in some other important cases such a decomposition does exist.

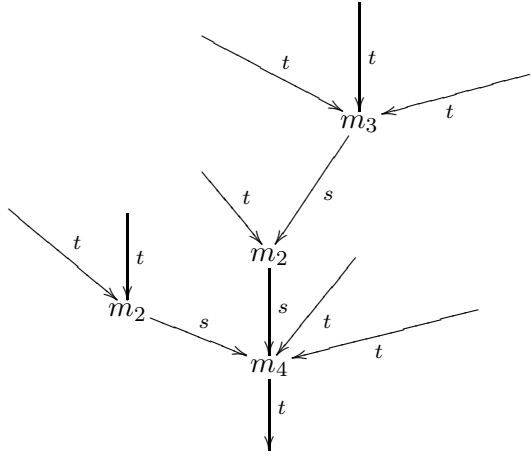
\begin{figure}[h]
\[
\xymatrix
@C=3.5ex@R=2.50ex@M=0.3EX
{
&&&&& \ar[ddd]^t \\
&&& \ar[ddrr]^t \\
&&&&&&&& \ar[dlll]^t\\
&&&&& m_3 \ar[dddl]^s\\
&&& \ar[ddr]^t \\
 \ar[ddrr]^t && \ar[dd]^t\\
 &&&& m_2 \ar[dd]^s & \ar[ddl]^t \\
 && m_2 \ar[drr]^s &&&&& \ar[dlll]^t \\
 &&&& m_4 \ar[dd]^t \\
 \\
 &&&& \\
}
\]
\caption{Definition of $\tilde{m}_T$.}
\label{planartree}
\end{figure}

We define new operators
$\tilde{m}_n:V ^{\otimes n}\rightarrow V, \quad n\geq 2$
as follows. Let $t:V \to V$ be the projection
onto $X$ along $\im(d)\oplus W$, and let $s:V \to V$ be inverse to $d$ on $W$ and $0$ on
 $X\oplus\im(d)$. Let $T$ be a planar rooted tree with $n+1$ extremities; we assume that each vertex has valence at least~$3$.
We label the extremities by $t$, all other edges by $s$ and each vertex of valence $v$ by $m_{v-1}$. Then
$\tilde{m}_T$ is constructed by working from the canopy of the tree down to the trunk, composing
labels in an obvious manner. For the tree pictured in Figure~\ref{planartree}
we have
\[
\tilde{m}_T=
tm_4
(\id \otimes sm_2 \otimes \id^{\otimes 2})
(sm_2\otimes \id\otimes sm_3 \otimes \id^{\otimes 2})
   t^{\otimes 8}.
\]

We define
$\tilde{m}_n = \sum_T  \tilde{m}_T: V^{\otimes n}\to V$, where
the sum is taken over all planar rooted trees with $n+1$ extremities.
It can be shown by direct calculation \cite{Mer, markl, KS} that
the $\tilde{m}_n$ satisfy the higher associativity conditions (\ref{eq:ainf}).
The compatibility (\ref{eq:cyclic}) of the $\tilde{m}_n$ with the bilinear form is
also known (by Lazaroiu \cite{Laz} in general and Kajiura \cite{Kajiura} when the form is non-degenerate). By restriction we obtain a cyclic $A_\infty$-structure on
$H(V)$, the minimal model, and it is independent (in an appropriate
sense) of the choice of decomposition of $V$.

Our goal is to give a conceptual explanation of this, extending our earlier paper \cite{CL2}.
So our results extend to algebras over arbitrary cyclic operads.

\section{Abstract Hodge decomposition}

\subsection{Definition}
Let $V$ be a dg vector space equipped with a bilinear form
$$\langle ,\rangle:V\times V \to \ground,$$
In addition, we assume that $\langle ,\rangle$ is homogeneous (i.e. either even or odd) and
graded symmetric or antisymmetric:
\begin{equation}\label{eq:formsymmetric}
\langle x,y\rangle = \pm(-1)^{|x||y|}\langle y,x\rangle,
\end{equation}
for any homogeneous $x,y\in V$.
The form is assumed to be compatible with the differential $d$ on $V$:
\begin{equation}\label{eq:dform}
\langle dx, y \rangle +(-1)^{|x|} \langle x, dy \rangle =0.
\end{equation}

We say that  the form $\langle,\rangle$ is \emph{non-degenerate} if for any nonzero $x\in V$ there exists
$y\in V$ such that $\langle x,y \rangle\neq 0$ and \emph{strongly non-degenerate} if it determines an isomorphism $V\rightarrow V^*$. Clearly a form can be strongly non-degenerate if and only if $V$ is finite dimensional and in this case this notion is equivalent to that of non-degeneracy.

For any dg subspace $W\subseteq V$, we denote by $W^\perp$ the dg subspace consisting of vectors $x$
\emph{orthogonal} to $W$, i.e. such that $\langle x,y\rangle=0$ for all $y\in W$. We say that
$W$ is \emph{isotropic} if $W\subseteq W^\perp$.

\begin{defi}\label{def_H}
A Hodge decomposition of $V$ is a pair of operators
$$t:V\to V \quad \text{and} \quad s:V\to V,$$
of even and odd degrees, respectively,
 such that
\begin{enumerate}
\item $s^2=0$,
\item $\langle sx, y \rangle = (-1)^{|x|} \langle x, sy \rangle$,
\item $sd+ds=1-t$,
\item $dt=td$,
\item $\langle tx, y \rangle = \langle x, ty \rangle$,
\item $t^2=t$,
\item $st=ts=0$,

\end{enumerate}
The Hodge decomposition is called \emph{trivial} if $t=\id_V$.
Following the suggestion of Dan Grayson we will call a Hodge decomposition \emph{harmonious} if $dt=0$; it was called `canonical' in \cite{CL2} but we feel that the latter terminology may be misleading.
\end{defi}

\begin{rem} \begin{itemize} \
\item The notion of a Hodge decomposition (without a biliniar form) is also known as \emph{splitting homotopy} or \emph{strong deformation retract} (SDR) data studied in the context of homological perturbation theory cf. \cite{BL}. It is clear that there is a direct sum decomposition $V=\Ker t\oplus\im t$ into orthogonal sub dg vector spaces; moreover $s$ determines a contractible homotopy on $\Ker t$.  The Hodge decomposition is then trivial if and only if $\im t=V$, and harmonious if and only if $\im t$ has zero differential and thus carries the homology of $V$.

\item In  view of (3) the operator $t$ is completely determined by $s$. It is easy to check that the conditions (4) and (5) are consequences of (1), (2), (3) and that (in the presence of (1), (2) and (3)) the conditions (6) and (7) together are equivalent to requiring that $sds=s$. Thus, a Hodge decomposition on $V$ is encoded in an odd self-adjoint operator $s$ of square zero and such that $sds=s$. It is harmonious if and only if $dsd=d$.
\end{itemize}
\end{rem}
\begin{defi}
We define an \emph{almost} Hodge decomposition in the same way as a Hodge decomposition, but just require that $ds+sd+t$ be invertible.
\end{defi}

It is always possible to `correct' an almost Hodge decomposition to a bona-fide one, as follows. The operator $ds+sd$ restricts to an automorphism of $\Ker t$, by assumption. Define the `Green operator' $G\in\End{V}$ to be inverse to $ds+sd$ on $\Ker t$ and $0$ on $\im{t}$. Since $d(ds+sd)=dsd=(ds+sd)d$ we see that $G^{-1}|_{{\Ker t}}$ commutes with $d$ and therefore so does $G$. Then we obtain
\[(sG)d+d(sG)=1-t,\]
i.e. a genuine Hodge decomposition with the operator $sG$ instead of $s$.
We remark that one automatically has $Gt=0$ and $Gs=sG$.
\begin{rem}\
\begin{itemize}\item
The significance of an almost Hodge decomposition is that the classical geometric Hodge decomposition of differential forms on a smooth or a K\"ahler manifold is an almost Hodge decomposition in our sense as will be explained below.
\item
Another potentially useful weakening of the notion of a Hodge decomposition is where only
conditions (2)-(5) are satisfied, except that the $t$ in (3) is replaced by $t^2$. Here
notice that $t$ is not determined by $s$, so (4) and (5) are not automatic consequences of (2) and (3). This version of a Hodge decomposition was employed in \cite{CL2}.
\end{itemize}
\end{rem}
The natural question is whether a Hodge decomposition of a dg vector space with a bilinear form always exists.
We first give the following preliminary result reformulating the notion of a Hodge decomposition in geometric terms.
\begin{prop}\label{geom}
Let $V$ be a dg vector space with an (anti)-symmetric bilinear form $\langle, \rangle$ as above. A harmonious Hodge decomposition on $V$ is equivalent to a decomposition of $V$ into a direct sum of three subspaces:
\[
V=\im d\oplus U\oplus W
\]
where $W\bot(\im d\oplus U)$, $W$ carries the homology of $V:H(V)\cong W$ and $U$ is an isotropic subspace of $V$.
\end{prop}
\begin{proof}
If the dg space $V$ is supplied with the operators $s$ and $t$ as in Definition \ref{def_H} we set $U=\im s$ and $W=\im t$; the required conditions are then checked in a straightforward fashion. Conversely, suppose that we are given a decomposition of $V$ in a direct sum of three spaces $\im d$, $U$ and $W$ as above. Since the space $\im d\oplus U$ is contractible it follows that $d$ maps $U$ isomorphically onto $\im d$ and we define $s|_{\im d}$ to be the inverse of $d:U\rightarrow \im d$ and $s|_{U\oplus W}=0$. Further, define $t$ to be the
projection onto $W$. It is immediate to check that $s^2=0$ and $sd+ds=1-t$. Furthermore, we have:
\begin{align*}\langle sv_1, v_2\rangle=&\langle sv_1,(sd+ds+t)v_2\rangle\\
=&\langle sv_1,dsv_2\rangle\\
=&(-1)^{|v_1|}\langle dsv_1, sv_2\rangle\\
=&(-1)^{|v_1|}\langle(1-sd-t)v_1,sv_2\rangle\\
=&(-1)^{|v_1|}\langle v_2, sv_2\rangle.
\end{align*} The remaining axioms for a Hodge decomposition are automatic.
\end{proof}
\begin{rem}
There is also an obvious and well-known analogue of the previous result without the presence of a form (which is actually a special case of it with the form being identically zero). In this situation a Hodge decomposition is equivalent to a direct sum decomposition $V=\im d\oplus U\oplus H(V)$ without any additional assumptions.
\end{rem}

\begin{theorem}\label{theorem:hodge}
Let $V$ be a dg vector space with an (anti)symmetric bilinear form $\langle, \rangle$ whose kernel has a dg complement and has finite codimension. Then a harmonious Hodge decomposition always exists.
\end{theorem}

\begin{proof}
Let $V^{\prime}$ be the kernel of $\langle,\rangle$, and choose a dg complement: $V=V^\prime\oplus V^{\prime\prime}$. Since
$V^{\prime}$ and $V^{\prime\prime}$ are orthogonal, it suffices to find harmonious Hodge decompositions for $V^{\prime}$
and $V^{\prime\prime}$ individually. Therefore we may assume that the form on $V$ is either identically zero or non-degenerate.

If the form is zero, then we first choose a subspace $W$ in $V$ representing $H(V)$ and then choose a complement $U$ to $\im d\oplus W$ in $V$; by the previous result this will give a desired Hodge decomposition.

Now suppose that the form is strongly non-degenerate (and thus, $V$ is finite dimensional). The first step is the same as before: we choose a subspace $W$ representing $H(V)$. Since our form is non-degenerate on $V$ it also non-degenerate on $W$ and therefore $V=W\oplus W^\bot$. Moreover, $\im d\in W^\bot$; it thus remains to find an isotropic complement to $\im d$ in $W^\bot$.

To this end, note that $\langle,\rangle$ can be viewed as an element in $(W^\bot\otimes W^\bot)^*$. The latter is a contractible dg vector space and so any cycle in it is in fact a boundary.  It follows that there exists another scalar product $(,)$ on $W^\bot$ such that $\langle v_1,v_2\rangle=(dv_1,v_2)+(-1)^{|v_1|}(v_1,dv_2)$ for any homogeneous $v_1,v_2\in W^\bot$. It follows that $(,)$ is non-degenerate on $\im d\in W^\bot$. Indeed, if there exists $v\in \im d$ such that for any $u\in W^\bot$ we have $(v,du)=0$, then $\langle v, u\rangle=0$ for any $u\in W^\bot$ contradicting the non-degeneracy of $\langle,\rangle$.

Set $U$ be the subspace of $W^\bot$ orthogonal to $\im d$ \emph{with respect to $(,)$}. Then, since $\im d\subset W^\bot$ is \emph{finite dimensional}, we have $W^\bot=\im d\oplus U$. It is also clear that $U$ is isotropic (with respect to the bilinear form $\langle, \rangle$) and this finishes the proof.
\end{proof}
\begin{rem}\
\begin{itemize}\item
In the above proof we used the finite-dimensionality assumption in two places: \begin{enumerate}\item
in the fact that $V=W^\bot\oplus W$ (this uses finite-dimensionality of $W=H(V)$) and
\item in the fact that $W^\bot=\im d\oplus U$ (this uses finite-dimensionality of $\im d$ and thus -- that of $W^\bot$)\end{enumerate}
In the most interesting cases provided by the classical geometric Hodge decomposition (see below) the space $H(V)$ is finite dimensional but $\im d$ is not and additional analytical arguments are needed in order to prove the existence of the desired complement. It is likely that there exist examples of maximal isotropic subspaces (such as $\im d$) in general infinite-dimensional spaces not having an isotropic complement; however we have been unable to construct such an example.
\item
In this paper we work with ${\mathbb Z}/2$-vector spaces; however all our results carry over almost verbatim to $\mathbb Z$-graded vector spaces. In this context, assuming that the given form $\langle, \rangle$ is \emph{homogeneous} and that the dg-vector space $V$ has finite-dimensional graded components Theorem \ref{theorem:hodge} still holds.
\end{itemize}
\end{rem}
\begin{example}\label{Hodge}\
\begin{enumerate}
\item Let $A$ be the algebra of differential forms on a smooth oriented compact manifold $M$. Then $A$ has a bilinear form: for two forms $\omega, \eta\in A$ set $\langle \omega, \eta\rangle=\int_M\omega\wedge \eta$; this form is clearly compatible with the external differential $d$. The form has zero kernel; however since $A$ is infinite dimensional Theorem \ref{theorem:hodge} is not applicable. The desired result follows instead from geometric Hodge decomposition, see e.g. \cite{Mor}. One introduces a Riemannian metric $g(,)$ on $M$, the associated volume form $d\mu$ and the operator $*$ determined by the formula $\omega\wedge *\eta=g(\omega,\eta)d\mu$. The operator $*$ has the following properties:
    \begin{equation}\label{ast}\omega\wedge *\eta=\eta\wedge *\omega;\end{equation}
     \begin{equation}\label{square}
     **\omega=(-1)^{|\omega|(n-|\omega|)}\omega.
     \end{equation}
     Next, there is defined a global scalar product $(,)$ on $A$ by the formula $(\omega,\eta)=\int_M\omega\wedge*\eta$ and an operator $d^*$ defined by $d^*(\omega)=(-1)^{n(|\omega|+1)+1}*d*$. A simple manipulation with formula (\ref{ast}) shows that
    \begin{equation}\label{adjoint}\langle d^*\omega,\eta\rangle=(-1)^{|\omega|}\langle \omega,d^*\eta\rangle.\end{equation} Then Hodge's theorem asserts that there is a direct sum decomposition
    \[A=H(A)\oplus \im d^*\oplus \im d\] where $H(A)$ is represented as $\Ker d\bigcap\Ker d^*$.
    It follows from integration by parts and formula (\ref{adjoint}) that $H(A)$ is orthogonal to
    $\im d^*\oplus \im d$ in the sense of the bilinear product $\langle,\rangle$. It further follows from formula (\ref{adjoint}) that the subspace $\im d^*$ is isotropic. Therefore by Proposition \ref{geom} there exists an abstract Hodge decomposition on $A$. Note that the Laplacian $dd^*+d^*d$ is not the identity on $\im d^*\oplus \im d$ but is an invertible operator (whose inverse is the Green operator), therefore the decomposition obtained most naturally is, in fact, what we termed an \emph{almost Hodge decomposition}. Of course, an almost Hodge decomposition is equivalent to a genuine one
\item Let $A=\Omega^{0,*}$ be the Dolbeault algebra of a complex Calabi-Yau manifold $M$ of dimension $n$. Thus, $A$ consists of differential forms of type $(0,i), i=0,\ldots,n$, that is those forms which could be written locally as $\sum f(z_{i_1},\ldots, z_{i_n},\bar{z}_{i_1},\ldots, \bar{z}_{i_n})dz_{i_1},\ldots, dz_{i_n}d\bar{z}_{i_1},\ldots, d\bar{z}_{i_n}$. The holomorphic volume form $\gamma\in \Omega^{n,0}$ determines a bilinear pairing $\langle\omega, \eta\rangle=\int_M\gamma\omega\eta$ compatible with the differential $\bar{\partial}$ in $A$. Indeed, we have for $\omega, \eta\in A$:
\[0=\int_M\bar{\partial}(\gamma\omega\eta)=\int_M\gamma\bar{\partial}(\omega\eta)=\langle\bar{\partial}(\omega),\eta\rangle+
(-1)^{|\omega|}\langle\omega,\bar{\partial}(\eta)\rangle.\]

As in the previous example this bilinear form is non-degenerate but $A$ is infinite dimensional. However  an abstract Hodge decomposition on $A$ exists and is given by the complex-analytic Hodge decomposition, see, e.g. \cite{GH}. Just as in the $C^\infty$-case one introduces
the operator $\bar{\partial}^*$ and obtains the following decomposition:
\[A=H(A)\oplus \im \bar{\partial}\oplus\im\bar{\partial}^* \]
which satisfies the conditions of Proposition \ref{geom} and so there results an abstract Hodge decomposition of $A$. Moreover, if $E$ is a holomorphic vector bundle (or a dg holomorphic vector bundle) then the dg-algebra of $\End(E)$-valued $(0,*)$-forms $A\otimes_{C^{\infty}(M)}\Gamma(\End(E))$ has a compatible bilinear pairing: for $f,g\in\Gamma(\End(E))$ and $\omega,\eta\in A$ we have
\[\langle\omega\otimes f, \eta\otimes g\rangle=\int_M\gamma\omega\eta \Tr(fg).\]
Again, geometric Hodge decomposition translates into an abstract Hodge decomposition on this algebra.
\end{enumerate}
\end{example}

\section{Cyclic operads and their algebras}
The notion of a \emph{cyclic operad} was introduced by Getzler and Kapranov \cite{GeK}, cf. also \cite{MSS}, Section 5.1.
\begin{defi}\ref{cyc}
An operad $\O$ is cyclic if the right $S_{n}$-action on $\O(n)$ extends to a right $S_{n+1}$-action such that:
\begin{enumerate}
\item
${\bf 1}^{z(1)}=\bf 1$ where $\bf 1\in\O(1)$ is the operadic unit;
\item for any $\phi\in\O(m)$ and $\psi\in\O(n)$
\begin{enumerate}
\item
$(\phi\circ_1\psi)^{z(m+n-1)}= \psi^{z(n)}\circ_n\phi^{z(m)};$
\item
$(\phi\circ_i\psi)^{z(m+n-1)}=\phi^{z(m)}\circ_{i-1}\psi ~{\text {for} } ~i=2,\ldots,m$
\end{enumerate}
\end{enumerate}
\end{defi}
We will also have a chance to use the notion of an \emph{anticyclic} operad; its definition is obtained from the one above by multiplying the right hand sides of formulas (1) and (2) by $-1$. Anticyclic operads appear as cobar-constructions of cyclic operads and this is the main reason for considering them.

Most known operads are cyclic. We now discuss the notion of an algebra over a cyclic operad. It is usually assumed that the underlying space of such an algebra possesses a strongly non-degenerate  inner product. In cases when such a product (or coproduct) is not strongly non-degenerate there are several inequivalent notions of a cyclic algebra.
\subsection{Different notions of algebras over cyclic operads}
Given a $dg$ space $V$, we denote by $E(V)$ the endomorphism operad on $V$, so that $E(V)(n)=\Hom(V^{\otimes n}, V)$. If $V$ is moreover equipped with a  symmetric bilinear form, we define the endomorphism cyclic operad $\E(V)$ by $\E(V)(n)=V^{\otimes n+1}$.

\begin{defi}\label{cyc}
Let $\O$ be a cyclic operad. Let $V$ be an dg space (possibly
infinite dimensional).
\begin{enumerate}
\item A \emph{standard} $\O$-algebra structure on $V$ is a symmetric bilinear form on $V$ together with a map of cyclic operads $\O\to\E(V)$.
\item A \emph{cyclic} $\O$-algebra structure on $V$ is a
symmetric bilinear form on $V$ together with a map of operads $\O\to E(V)$ inducing $S_{n+1}$-module maps
$\O(n)\to (V^{\otimes n+1})^*$. We will also say that the given form on  $V$ is \emph{invariant} with respect to the $\O$-algebra structure.
\item A \emph{Casimir} $\O$-algebra structure on $V$ is a tensor in $V\otimes V$ together with a map of operads $\O\to E(V)$ inducing $S_{n+1}$-invariant maps $\O(n)\to V^{\otimes n+1}$.
\end{enumerate}
\end{defi}

\begin{rem}
If the form on $V$ is strongly non-degenerate, then the notions of standard and cyclic algebras are equivalent. Furthermore, if the tensor in the definition of Casimir algebra induces an isomorphism $V^*\cong V$, then the notion of Casimir algebra also reduces to either of the other two. In this situation $E(V)$ is a cyclic operad and a cyclic $\O$-algebra structure on $V$ is simply a map of cyclic operads $\O\rightarrow E(V)$.

In general the three notions are distinct. Standard algebras are the most natural from a purely abstract point of view, in the setting of cyclic and modular operads. However it is unclear whether interesting examples of standard algebras with degenerate forms are plentiful in nature. Note that a standard algebra over the associative cyclic operad is a dg vector space $V$ together with a bilinear form and a three-tensor $\sum a_i\otimes b_i\otimes c_i\in V^{\otimes 3}$, subject to a quadratic relation; it is thus not even an associative algebra.

The notion of Casimir $\O$-algebra is introduced in Hinich-Vaintrob \cite{HV}.

Dual to the Casimir algebras are the cyclic algebras, in which we are most interested. The cyclic algebras are the operadic algebras equipped with invariant bilinear form, as defined by Getzler and Kapranov \cite{GeK}.
 Note that the induced maps $\O(n)\to (V^*)^{\otimes n+1}$ automatically intertwine the action of $S_n$. Hence it is enough to further ensure that the action of the cyclic group
 $\langle z(n)\rangle$ is preserved.
\end{rem}

Let $\O$ be a cyclic operad and $V$ be an $\O$-algebra (viewing $\O$ as just an operad). Then a bilinear form $\langle,\rangle$ on $V$ determines a cyclic $\O$-algebra structure on $V$ if and only if the map $B_n:\O(n)\otimes V^{\otimes(n+1)}\rightarrow \ground$ defined by the formula
\[B_n(\phi\otimes v_0 \otimes \ldots \otimes v_{n})
= \langle \phi(v_0, \ldots , v_{n-1}), v_{n}\rangle\]
is invariant under the action of $S_{n+1}$ on $\O(n) \otimes V^{\otimes(n+1)}$.
Note that the symmetry of the form is a consequence of the definition, if the operad is unital.

Our main result turns on the following simple but important lemma:

\begin{lem}\label{lem_cyclic}
Suppose that a cyclic operad $\O$ is generated by a sub- $\S_+$-module $M$.
Then a bilinear form $\langle,\rangle : V\times V \to \ground$
on an $\O$-algebra $V$ is invariant if and only if
for all $n\geq 0$ the restriction of $B_n$ to
$M(n)\otimes V^{\otimes(n+1)}$ is invariant under the action of $S_{n+1}$.
\end{lem}

\begin{proof}
We need to prove that for all $\phi\in\O(n)$ and all
$v_0,\ldots,v_{n}\in V$
 we have
$$\langle \phi(v_0,\ldots,v_{n-1}), v_{n}\rangle
=(-1)^{|v_0|(|v_1|+\ldots+|v_{n}|)}
\langle \phi^{z(n)}(v_1,\ldots,v_{n}), v_0 \rangle.$$

So it suffices to show that if this property holds for $\phi\in\O(n)$ and
$\psi\in\O(m)$ then it is true for a composition $\psi\circ_i\phi\in\O(m+n-1)$. For simplicity we assume that the elements $v_0,\ldots,v_n$ are even; the proof in the general case is slightly messier because of the signs. First suppose $i>1$. Then
\begin{eqnarray*}
\langle  (\psi\circ_i\phi)(v_0,\ldots,v_{m+n-2}), v_{m+n-1} \rangle
& = &
\langle  \psi(v_0,\ldots,v_{i-2}, \phi(v_{i-1},\ldots,v_{i+n-2}),v_{i+n-1},\ldots,v_{m+n-2}), v_{m+n-1} \rangle \\
& = &
\langle \psi^{z(m)}(v_1,\ldots,v_{i-1}, \phi(v_i,\ldots,v_{i+n-1}),v_{i+n},\ldots,v_{m+n-1}) , v_0 \rangle \\
&=&
\langle (\psi^{z(m)}\circ_{i-1}\phi)(v_1,\ldots,v_{m+n-1}), v_0 \rangle \\
& = &
\langle  (\psi\circ_i\phi)^{z(m+n-1)}(v_1,\ldots,v_{m+n-1}), v_0 \rangle,
\end{eqnarray*}
where in the second line we used the invariance  for $\psi$ while the last line used the cyclicity of the operad $\O$ (property (2b) of Definition \ref{cyc}).

For $i=1$ we have:
\begin{eqnarray*}
\langle  (\psi\circ_1\phi)(v_0,\ldots,v_{m+n-2}), v_{m+n-1} \rangle
& = &
\langle \psi(\phi(v_0,\ldots,v_{n-1}),v_{n},\ldots,v_{m+n-2}), v_{m+n-1} \rangle \\
& = &
\langle \psi^{z(m)}(v_{n},\ldots,v_{m+n-1}), \phi(v_0,\ldots,v_{n-1})  \rangle \\
& = &
\langle \phi(v_0,\ldots,v_{n-1}),
\psi^{z(m)}(v_{n},\ldots,v_{m+n-1}) \rangle \\
& = &
\langle \phi^{z(n)}(v_1,\ldots,v_{n-1}, \psi^{z(m)}(v_{n},\ldots,v_{m+n-1})),
v_0 \rangle \\
& = &
\langle (\phi^{z(n)}\circ_n\psi^{z(m)})(v_1,\ldots,v_{m+n-1}), v_0 \rangle \\
& = &
\langle (\psi\circ_1\phi)^{z(m+n-1)}(v_0,\ldots,v_{m+n-1}), v_0 \rangle. \\
\end{eqnarray*}

Here the second line used the invariance of $\psi$, the fourth -- the invariance of $\phi$ and the last -- the cyclicity of $\O$ (property (2a) of Definition \ref{cyc}).
\end{proof}

\section{Main theorem: an explicit form of a minimal model of cyclic algebras}
In this section we assume that a cyclic operad $\O$ is \emph{admissible}
i.e. that $\O(n)$ is finite dimensional for all $n$ and that $\O(1)$ is a one-dimensional space spanned by the operadic unit $\bf 1$. The cobar-construction $\B\O$ of $\O$ (or, more pedantically, that of the dual cooperad $\O(n)^*$) is an anticyclic operad, see \cite{GeK}.

Our main result relies heavily on our previous work \cite{CL2} and we now recall some of the constructions and terminology of this paper. Associated to $\O$ is its canonical cofibrant resolution (essentially the double cobar-construction) $\bv\O$ and another operad $\BV\O$ whose algebras are $\O$-algebras together with a choice of a Hodge decomposition (no invariant inner product is assumed yet). For an $\O$-algebra structure on a dg vector space $V$ we have the following commutative diagram:
\begin{equation}\label{operadmaps}\xymatrix {\bv\O\ar^i[r]\ar_p[d]&\BV\O\ar@{-->}^h[d]\ar^q@/^/[dl]\\ \O\ar_f[r]\ar^j@/^/[ur]&E(V)& E(H(V))\ar@{_{(}->}[l]}\end{equation}
Here $\BV\O$ is given by freely adjoining to $\O$ an odd element $s$ and an even element $t$, both of arity $1$ and subject to the relations $0=s^2=st=ts=t^2-t$ and with differential extending that on $\O$ and such that $d(s)={\bf 1}-t$ and $d(t)=0$. The operad $\BV\O$ can be viewed as the operad of $\O$-decorated trees whose internal edges are marked by the symbols $s$ or $t$; the vertices, except those adjacent to the extremities are required to be at least trivalent. There results an inclusion $j:\O\rightarrow \BV\O$ obtained by viewing $\O$ as the space of trees with no internal edges; this inclusion is split by the projection $q:\BV\O\rightarrow \O$ defined by setting $s\mapsto 0,t\mapsto \bf 1$.

The operad $\bv\O$ is the suboperad in $\BV\O$ consisting of $\O$-decorated trees whose extremities are connected to $t$-edges by bivalent vertices, and $p=q\circ i$. The map $f$ is the given $\O$-algebra structure on $V$ and the dotted arrow $h$ corresponds to a choice of a Hodge decomposition on $V$. The inclusion $E(H(V))\hookrightarrow E(V)$ also comes from a Hodge decomposition on $V$ (which we always assume to be harmonious).

The action of $\bv\O$ restricts to the image $\im t$ of the operator $t$ of the Hodge decomposition thus giving a $\bv\O$-algebra structure on $H(V)$.
In the case when the map $\bv\O\rightarrow\O$ has a splitting, as is the case when $\O$ itself is a (suspended) cobar-construction, this $\bv\O$-algebra structure pulls back to an $\O$-algebra structure on $H(V)$  giving a \emph{minimal model} of the $\O$-algebra $V$. The structure maps of this canonical (modulo the choice of a Hodge decomposition) minimal model are given by a Merkulov-like formula, see \cite{CL2}, Theorem 3.2.

Moreover, any two choices of a Hodge decomposition give \emph{homotopy equivalent} minimal models by Theorem 4.14 of \cite{CL2}. Let us recall the notion of homotopy equivalence of operadic algebras. It is based on a special dg algebra $D=\ground [z,dz]$ of polynomial forms on the interval. Note that there are two restriction maps $\ev_{0,1}:D\rightarrow \ground$ corresponding to the restrictions of differential forms to $0$ and $1$. Associated to any (cyclic) operad $\O$ is a (cyclic) operad $\O\otimes D$ given by $\O\otimes D(n)=\O(n)\otimes D$.
\begin{defi}\
\begin{enumerate}
\item
Two maps between operads $f_0,f_1:\O\rightarrow\O^\prime$ are \emph{homotopic} if there exists a map $f:\O\rightarrow \O^\prime\otimes D$ such that $f_0=(id\otimes \ev_0)\circ f$ and $f_1=(id\otimes \ev_1)\circ f$. If $\O$ and $\O^\prime$ are cyclic operads then $f_0, f_1$ are \emph{cyclically homotopic} if $f$ a map of cyclic operads.
\item
Let $V$ be a dg vector space supplied with a bilinear form giving a map of $S$-modules $E(V)\to\E(V)$. Two $\O$-algebras corresponding to operad maps $f_0:\O\rightarrow E(V)$ and $f_1:\O\rightarrow E(V)$ are called (cyclically) \emph{homotopy equivalent} if $f_0$ and $f_1$ are homotopic via a homotopy $\O\to E(V)\otimes D$ and the induced maps $\O(n)\to E(V)\otimes D\to\E(V)\otimes D$ are maps of $S_{n+1}$-modules.
\end{enumerate}
\end{defi}

The operad $\BV\O$ has an obvious cyclic structure with generators $s$ and $t$ being invariant with respect to the action of $S_2$; the suboperad $\bv\O$ is then also a cyclic operad and the maps $i,j,p,q$ clearly are maps of cyclic operads.

Now suppose that $V$ has a symmetric bilinear form, not necessarily non-degenerate, making it into a cyclic $\O$-algebra. The following result
asserts that Hodge decompositions allow one to build minimal models in a unique up to homotopy way.
\begin{theorem}\label{main}Let $V$ be a cyclic $\O$-algebra supplied with a Hodge decomposition.
\begin{itemize}
\item The
induced $\bv\O$-algebra structure on $tV$ makes it a cyclic $\bv\O$-algebra.
The cyclic $\bv\O$-algebra structures on $H(V)$ induced by any two harmonious Hodge decompositions are cyclically homotopy equivalent.
\item
Suppose that $\O=\B\P$, the cobar-construction of an anticyclic operad $\P$. Then the Merkulov-type formula as in Theorem 3.2 of \cite{CL2} gives a minimal model for $V$ as a cyclic $\O$-algebra.
Moreover, the cyclic $\O$-algebra structures on $H(V)$ induced by any two harmonious Hodge decompositions are cyclically homotopy equivalent.
\end{itemize}
\end{theorem}

\begin{proof}A choice of Hodge decomposition on $V$ allows one to extend the $\O$-algebra structure on $V$ to a $\BV\O$-algebra structure. Now $\BV\O$ is generated by the sub-$\S_+$-module $K=\O\oplus \ground s \oplus \ground t$, and the induced maps $K(n)\to(V^{\otimes n})^*$ are $S_{n+1}$-module homomorphisms as long as the Hodge decomposition is chosen compatible with the bilinear form. So by Lemma \ref{lem_cyclic} the  $\BV\O$-algebra structure on $V$ is cyclic. It follows that the
$\bv\O$-algebra structure obtained by restriction is also cyclic.

The homotopy between different minimal models described in \cite{CL2} (Theorem 4.14) is easily checked to be a cyclic map: the image of $s$ and $t$ are compatible with the extension of the inner product to the larger coefficient dg agebra $D=\ground[z,dz]$.

The second statement about minimal models of $V$ as an $\O$-algebra follows from the first since the canonical map
$\O=\B\P\rightarrow \bv\O=\B\B\B\P$ (given by applying the contravariant functor $B$ to the canonical resolution $\B\B\P\rightarrow \P$) is a map of cyclic operads.
\end{proof}
There is an important special case where a cyclic algebra $V$ over $\O$ has finite-dimensional homology: $\dim H(V)<\infty$ and the given bilinear form on $H(V)$ is non-degenerate.
In this situation the endomorphism operad $E(H(V))$ is cyclic and a minimal model of $V$ is specified by a map of \emph{cyclic} operads $\bv\O\rightarrow E(H(V))$. There results a map of \emph{modular operads} (cf. \cite{GeK1} concerning this notion): $\overline{\bv\O}\rightarrow E_{mod}(H(V))$ where the bar indicates the modular closure, i.e. the modular operad freely generated by the corresponding cyclic operad and $E_{mod}(H(V))$ denotes the modular endomorphism operad i.e.
$E_{mod}(H(V))((g,n))=\Hom(H(V)^{\otimes n-1}, H(V))$. We also get  a homotopy $\overline{\bv\O}\rightarrow E_{mod}(H(V))\otimes D$ corresponding to two choices of a Hodge decomposition, giving rise to a homotopy between modular operads.
We obtain the following corollary:
\begin{cor}
Let $\O$ and $V$ be as in Theorem \ref{main} and, in addition, suppose $H(V)$ is finite dimensional and the given bilinear form is non-degenerate on $H(V)=tV$.
\begin{itemize}
\item
The space $H(V)$ has the structure of an algebra over the modular operad
$\overline{\bv\O}$.
The  $\overline{\bv\O}$-algebra structures on $H(V)$ induced by any two compatible harmonious Hodge decompositions are homotopy equivalent as modular algebras.
\item
Suppose that $\O=\B\P$, the cobar-construction of an anticyclic operad $\P$. Then the Merkulov-type formula as in Theorem 3.2 of \cite{CL2} gives a minimal model for $V$ as a modular $\overline{\O}$-algebra.
Moreover, the modular $\overline{\O}$-algebra structures on $H(V)$ induced by any two compatible harmonious Hodge decompositions are  homotopy equivalent.
\end{itemize}
\end{cor}
\noproof
\begin{rem}
The above corollary allows one to extend the definition of characteristic classes (Kontsevich's `direct construction' formulated originally for $A_\infty$-algebras, cf. \cite{K1}, \cite{HL}) for cobar construction cyclic operads. Let $\O=\B\P$ be the cobar construction of an anticyclic operad $\P$ and let $\overline{\O}$ be its modular closure. Let $V$ be a cyclic algebra over $\O$ such that $V$ possesses a Hodge decomposition and the induced form on $H(V)$ is strongly non-degenerate. While $V$ itself is not necessarily an algebra over $\overline{\O}$, its homology is, and any two minimal model structures are weakly equivalent by the theorem.

Thus we have a map of modular operads $\overline{\O}\ra E_{mod}(H(V)))$. The characteristic class is given by the restriction of this map to the vacuum part $\overline{\O}((0))$ and therefore determines a cohomology class $[V]:=[H(V)]\in \overline{\O}((0))^*$ is well-defined and does not depend on the choice of a Hodge decomposition.

We remark that there is problem in carrying the above argument through for an arbitrary (i.e. not of the form $\B\P$) cyclic operad $\O$. Even though $\overline{\BV\O}\cong \BV(\overline{\O})$, the inclusion $\overline{\bv\O}\hookrightarrow
\bv(\overline{\O})$ might not be a quasi-isomorphism (even on the vacuum). For degenerate algebras the characteristic class construction would take values only in $\overline{\bv\O}((0))^*$, rather than in $\bv(\overline{\O})((0))^*$.
\end{rem}
\begin{example}\
\begin{enumerate}
\item
Let $V$ be the Sullivan minimal model of the algebra of differential forms on a simply connected Poincare duality space $M$. We saw that $V$ has a Hodge decomposition compatible with the Poincar\'e pairing. The resulting cyclic $C_\infty$-model of $V$ is what was called in \cite{La} a \emph{Stasheff model} and derived using different methods.
If $M$ is in fact a smooth manifold then we can alternatively use the geometric Hodge decomposition of the de Rham algebra of $M$ to obtain the same result.
\item
Consider the bounded derived category $D^b(M)$ of a smooth projective Calabi-Yau manifold $M$. It has a dg model ${\mathcal C}(M)$ obtained by taking the Dolbeault resolutions of coherent sheaves. The Hochschild cohomology of this category is, essentially, $H^*(M)$, the Hodge cohomology of $M$. On the other hand, the existence of a strong generator of $D^b(M)$, cf. \cite{BV} (represented by a complex of holomorphic vector bundles on M) implies that there is a dg algebra $V$ whose derived category is equivalent to $D^b(M)$.
As explained in Example \ref{Hodge} (3) Serre duality gives an invariant scalar product on $V$ and geometric Hodge decomposition allows one to construct a minimal cyclic model $(H(V), \{m_i\})$ for $V$. The Hochschild cohomology of $(H(V), \{m_i\})$ is isomorphic to the Hochschild cohomology of $V$ and thus -- to $H^*(M)$.

Recall the result \cite{Cos, KS1}
which states that the Hochschild cohomology of a cyclic $A_\infty$-algebra with a non-degenerate form supports an action of the chain operad of moduli spaces of Riemann surfaces.
 We conclude that the $H^*(M)$ forms an algebra over the chain operad of moduli of Riemann surfaces.
\end{enumerate}
\end{example}

\end{document}